\documentclass[sn-mathphys-num]{sn-jnl}

\usepackage{graphicx}%
\usepackage{multirow}%
\usepackage{amsmath,amssymb,amsfonts}%
\usepackage{amsthm}%
\usepackage{mathrsfs}%
\usepackage[title]{appendix}%
\usepackage{xcolor}%
\usepackage{textcomp}%
\usepackage{manyfoot}%
\usepackage{booktabs}%
\usepackage{algorithm}%
\usepackage{algorithmicx}%
\usepackage{algpseudocode}%
\usepackage{listings}%

\theoremstyle{thmstyleone}%
\newtheorem{theorem}{Theorem}
\theoremstyle{thmstyletwo}%

\theoremstyle{thmstylethree}%
\newtheorem{definition}{Definition}%
\newtheorem{conclusion}{Conclusion}%
\begin{document}

\title[Investigation of Determinants of Fibonacci-Hessenberg-Lorentz Matrices and Special Number Sequences]{Investigation of Determinants of Fibonacci-Hessenberg-Lorentz Matrices and Special Number Sequences}


\author*[1]{\fnm{Ibrahim} \sur{G\u{o}kcan}}\email{gokcan@artvin.edu.tr}

\author[2]{\fnm{Ali Hikmet} \sur{De\u{g}er}}\email{ahikmetd@ktu.edu.tr}
\equalcont{These authors contributed equally to this work.}

\affil*[1]{\orgdiv{Faculty of Arts and Sciences}, \orgname{Artvin \c{C}oruh University},\city{Artvin}, \postcode{08000},\country{T\"{u}rkiye}}

\affil[2]{\orgdiv{Departmant of Mathematics}, \orgname{Karadeniz Technical University},\city{Trabzon}, \postcode{61000},\country{T\"{u}rkiye}}


\abstract{The research aims to construct a new type of matrix called the Fibonacci-Hessenberg-Lorentz matrix by multiplying Fibonacci-Hessenberg matrices with Lorentz matrix multiplication. The study will start by examining the properties of Hessenberg and tridiagonal matrices and then focus on developing the Fibonacci-Hessenberg matrix using Fibonacci sequences. By multiplication it with a Lorentz matrix multiplication, the resulting matrix, the Fibonacci-Hessenberg-Lorentz matrix, will be analyzed to obtain special number sequences through its determinants for $n\geq1$. 

The primary objective is to explore whether the determinants of these matrices can generate new or known number sequences, where the elements are expressed as functions of the matrix parameters. Furthermore, the research will attempt to generalize these sequences of using Fibonacci numbers to establish a generalized formula for their terms. Ultimately, the goal is to derive a mathematical representation that connects the characteristics of the newly defined matrices to well-known special sequences in mathematics.
}




\maketitle

\section{Introduction}\label{sec1}
\subsection*{Fibonacci and Lucas Sequences}

Many sequences can be generated with different initial conditions and recurrence relations. The Fibonacci sequence, one of the most studied sequences in the literature, can broadly be defined by recurrence relation $F_{n}=F_{n-1}+F_{n-2}$ and initial conditions $F_1=1$ and $F_{2}=1$. That is, the terms of Fibonacci sequence can be given as $1,1,2,3,\cdots$. The approach of recurrence relation used in Fibonacci number sequence is similar to that used in Lucas number sequence. In the Lucas sequence, relation $L_{n}=L_{n-1}+L_{n-2}$ and initial conditions $L_{1}=2$ and $L_{2}=1$ are used and, sequence $2,1,3,4,\cdots$ is achieved. Numerous studies have attempted to investigate number sequences; especially [$\color{blue}{1-10}$] can be given as references. In [$\color{blue}{11}$], generalization of Lucas numbers is demonstrated and a connection is attained between   $k$-Lucas numbers and Fibonacci numbers. In addition, generalized $k$-Fibonacci and Lucas numbers are analyzed,  Binet formulas are investigated for these sequences in [$\color{blue}{12}$]. In [$\color{blue}{13}$], a closed formula are accessed for Horadam polynomials and polynomial sequences by using tridiagonal determinants. A closed formula giving the elements of Horadam sequence are achieved by tridiagonal matrices in [$\color{blue}{14}$]. The formulas related to $(p,q,r)$-tribonacci polynomials and generalized tribonacci sequences are presented by using Hessenberg matrices in [$\color{blue}{15}$]. 
\subsection*{Hessenberg and Tridiagonal Matrices and its Determinants}
The term $\alpha_{ij}$ is used here to refer to an element in a matrix. Several definitions of Hessenberg matrix have been proposed. First, Hessenberg matrix is defined as matrix $n\times n$ such that $\alpha_{j,j+1}\neq 0$ and $\alpha_{ij}=0$ for $j>i+1$ but not lower triangular. Second, main diagonal and its lower triangular matrix are compose of elements different from $0$, but upper triangular matrix compose of elements $0$. The following matrix is a good illustration of Hessenberg matrix:
\begin{equation}
H_{n}=\begin{pmatrix}
h_{1,1} & h_{1,2}& \cdots & 0 & 0 \\
a_{2,1} & a_{2,2}& \cdots & 0 & 0 \\
\vdots & \vdots& \ddots & 0 & 0 \\
h_{n-1,1} & h_{n-1,2}& \cdots & h_{n-1,n-1}& h_{n-1,n} \\
h_{n,1} & h_{n,2}& \cdots & h_{n,n-1}& h_{n,n} 
\end{pmatrix}
\end{equation}
Assume that the sequence consisting of the determinants of the Hessenberg matrix  $H_{n}$ is $\lbrace \vert H_{n}\vert,n\geqslant1 \rbrace$. For instance, for $n=1$, $H_{n}=h_{1,1}$. For $n\geqslant2$, the general formula giving determinants is:
\begin{equation}
\vert H_{n}\vert=h_{n,n}\vert H_{n-1}\vert+\sum_{r=1}^{n-1}\left[\left(-1\right)^{n-r}h_{n,r}\prod_{j=r}^{n-1}h_{j,j+1}\vert H_{r-1}\vert \right] 
\end{equation}
There are multiple definitions of a tridiagonal matrix. One of these can be given as if the elements of different from $0$ only on main diagonal, upper diagonal and lower diagonal, matrix is defined as tridiagonal matrix. Mathematically, the tridiagonal matrix can be defined as $\alpha_{ij}=0$ for $j>i+1$ and $i>j+1$,  $\alpha_{ii}\neq0$ and $\alpha_{j,j+1}\neq0$ for in other situations.\\
The following matrix illustrates a tridiagonal matrix clearly.
\begin{equation}
T_{n}=\begin{pmatrix}
t_{1,1} & t_{1,2}& \cdots & 0 & 0 \\
t_{2,1} & t_{2,2}& \cdots & 0 & 0 \\
\vdots & \vdots& \ddots & 0 & 0 \\
0 & 0 & \cdots & t_{n-1,n-1}& t_{n-1,n} \\
0 & 0 & \cdots & t_{n,n-1}& t_{n,n} 
\end{pmatrix}
\end{equation}
However, the set of Hessenberg matrices cover the set of tridiagonal matrices. Therefore, obtaining determinant of a tridiagonal matrix can be used above equation. For some special values, the determinants of the Hessenberg and Tridiagonal matrices were associated with the general terms of the Fibonacci number sequence. This finding broadly supports the work of other studies in this area linking Hessenberg and Tridiagonal matrices with Fibonacci number sequence. A considerable amount of literature has been published on Hessenberg matrices whose determinants are related to Fibonacci numbers. Here, [$\color{blue}{16-20}$] can be given as a reference to some studies.\\
Here, some matrix examples and findings that have been extensively examined in the literature are given. In [$\color{blue}{16}$], it is shown that the determinant of the matrix, whose elements on the main diagonal are $1$, elements on the upper and lower diagonals are $i$, and other elements are $0$, is $F_{n+1}$ for $n\in\mathbb{N}$.\\
This matrix has been shown as follows with matrix $A_{n}$:
\begin{equation}
A_{n}=\begin{pmatrix}
1 & i & \cdots & 0 & 0 \\
i & 1 & \cdots & 0 & 0 \\
\vdots & \vdots& \ddots & 0 & 0 \\
0 & 0 & \cdots & 1 & i \\
0 & 0 & \cdots & i & 1 
\end{pmatrix}
\end{equation}
In the present study, $b_{ij}$ is defined as an element in matrix $B_{n}$. In addition, assume that the elements of main diagonal are $2$ but $b_{nn}=1$, the elements of upper diagonal matrix and lower diagonal are $1$. Then, matrix $B_{n}$ can be demonstrated in the following:\\
\begin{equation}
B_{n}=\begin{pmatrix}
2 & 1 & \cdots & 0 & 0 \\
1 & 2 & \cdots & 0 & 0 \\
\vdots & \vdots& \ddots & 0 & 0 \\
1 & 1 & \cdots & 2 & 1 \\
1 & 1 & \cdots & 1 & 1 
\end{pmatrix}
\end{equation}
Matrix $B_{n}$ can be entitled upper Hessenberg matrix because of the elements of upper diagonal matrix are $0$ and main diagonal and elements of lower diagonal matrix are different from $0$. One interesting finding relevant to is the determinants of matrix $B_{n}$ are $\vert B_{1} \vert=1,\vert B_{1} \vert=1$ and $\vert B_{n} \vert=F_{n}$ for respectively $1,2$ and $n$.\\
In contrast to other studies, [$\color{blue}{16}$] in the determinants of Hessenberg matrices whose elements are complex numbers are associated with Fibonacci sequences. In [$\color{blue}{19}$], the authors examined the determinants of tridiagonal matrices whose elements are real numbers and associated determinants with Fibonacci sequences.\\
In previous studies, the first and second type Fibonacci numbers, $F_{n-1} t+F_{n-2}$ and $F_{n-1}+F_{n-2}t$  respectively, are found by giving a variable $t$ in the elements $\left(1,1\right)$ or $\left(n,n\right)$  of some Hessenberg and tridiagonal matrices.
\begin{definition}
Assume that $t$ is a real or complex number and $n$ is an integer. The terms the first and second type Fibonacci numbers, $F_{n-1} t+F_{n-2}$ and $F_{n-1}+F_{n-2}t$  respectively are used in its broadest sense to refer to all the determinant of the matrices created depending on $t$ and $n$. In the present study, the first and second type Fibonacci numbers is defined as $\left(t,n\right)$-Fibonacci.
\end{definition}
\begin{definition}
Assume that $H_{n}$ is a Hessenberg matrix and $t$ is an integer providing $m>0$ such that $\forall n \geqslant  m$. Hessenberg $H_{n}$ matrix is named as Fibonacci-Hessenberg matrix if determinants of Hessenberg $H_{n}$ matrix are $\left(t,n\right)$-Fibonacci numbers.
\end{definition}
Here, a well-known examples of some Fibonacci-Hessenberg matrices and sequences created by its determinants with $\left(t,n\right)$-Fibonacci numbers are presented for $n\geqslant1$.\\
First, a classic example of Hessenberg matrix is $C_{n,t}$ that its determinants are shown Fibonacci numbers:
\begin{equation}
C_{n,t}=\begin{pmatrix}
2 & 1 & \cdots & 0 & 0 \\
1 & 2 & \cdots & 0 & 0 \\
\vdots & \vdots& \ddots & 0 & 0 \\
1 & 1 & \cdots & 2 & 1 \\
1 & 1 & \cdots & 1 & t+1 
\end{pmatrix}
\end{equation}
For $n\geqslant1$, determinants of matrix $C_{n,t}$ are indicated as follows:\\
$\vert C_{1,t}\vert=t+1,\vert C_{2,t} \vert=2t+1,\vert C_{3,t}\vert=3t+2,\cdots,\vert  C_{n,t}\vert=F_{n+1} t+F_{n}$
$C_{n,t}$ is the first type $\left(t,n\right)$-Fibonacci number. Further analysis shows that the sequences $\lbrace 0,-1,-1,\cdots,-F_{n-1},\cdots\rbrace$, $\lbrace 1,1,2,\cdots,F_{n},\cdots\rbrace$,$\lbrace 2,3,5,\cdots,F_{n+2},\cdots\rbrace$ and $\lbrace 3,5,8,\cdots,F_{n+3},\cdots\rbrace$ are achieved for respectively $t=-1,0,1$ and $2$.\\
Second, matrix $D_{n,t}$ is obtained by replacing which element are $1$ in the upper diagonal of the matrix $C_{n,t}$ with $-1$. In the following, the matrix $D_{n,t}$ and its determinants are demonstrated: 
\begin{equation}
D_{n,t}=\begin{pmatrix}
2 & -1 & \cdots & 0 & 0 \\
1 & 2 & \cdots & 0 & 0 \\
\vdots & \vdots& \ddots & 0 & 0 \\
1 & 1 & \cdots & 2 & -1 \\
1 & 1 & \cdots & 1 & t+1 
\end{pmatrix}
\end{equation}
$\vert D_{1,t}\vert=t+1,\vert D_{2,t}\vert=2t+3,\vert D_{3,t}\vert=3t+5,\cdots,\vert D_{n,t}\vert=F_{2n-1} t+F_{2n}$
The results obtained from analysis of determinants of matrix $D_{n,t}$ are $\lbrace 0,1,2,\cdots,F_{2n-2},\cdots\rbrace$, $\lbrace 1,3,5,\cdots,F_{2n},\cdots,\cdots\rbrace$,$\lbrace 2,5,8,\cdots,F_{2n+1},\cdots,\cdots\rbrace$  and $\lbrace 3,7,11,\cdots,L_{2n},\cdots,\cdots\rbrace$ for respectively $t=-1,0,1$ and $2$.\\
Third, the other Fibonacci-Hessenberg matrix is matrix $E_{n,t}$ achieving by replacing which element is $2$ in the $\left(1,1\right)$ of the matrix $D_{n,t}$ with $1$. The matrix $E_{n,t}$ can be illustrated and generalized determinants can be presented as follows:
\begin{equation}
E_{n,t}=\begin{pmatrix}
1 & -1 & \cdots & 0 & 0 \\
1 & 2 & \cdots & 0 & 0 \\
\vdots & \vdots& \ddots & 0 & 0 \\
1 & 1 & \cdots & 2 & -1 \\
1 & 1 & \cdots & 1 & t+1 
\end{pmatrix}
\end{equation}
$\vert E_{1,t}\vert=t+1,\vert E_{2,t}\vert=2t+3,\vert E_{3,t}\vert=3t+5,\cdots,\vert E_{n,t}\vert=F_{2n-2} t+F_{2n-1}$
From operations of determinant we can see that sequences $\lbrace 0,1,2,\cdots,F_{2n-3},\cdots\rbrace$, $\lbrace 1,3,5,\cdots,F_{2n-1},\cdots\rbrace$, $\lbrace 2,5,8,\cdots,F_{2n},\cdots\rbrace$ and $\lbrace 3,7,11,\cdots,L_{2n-1},\cdots\rbrace$ are achieved for respectively $t=-1,0,1$ and $2$.\\
Fourth, Fibonacci-Hessenberg matrix $F_{n,t}$ is found by replacing which element is $1$ in the $\left(1,1\right)$ of the matrix $E_{n,t}$ with $2$. In addition, lower diagonal matrix are compose of elements of $-1$ and $1$ respectively. A matrix representing the matrix $F_{n,t}$ is given below:
\begin{equation}
F_{n,t}=\begin{pmatrix}
2 & -1 & \cdots & 0 & 0 \\
-1 & 2 & \cdots & 0 & 0 \\
\vdots & \vdots& \ddots & 0 & 0 \\
-1 & 1 & \cdots & 2 & -1 \\
1 & -1 & \cdots & -1 & t+1 
\end{pmatrix}
\end{equation}
For $n\geqslant1$, determinants of matrix $F_{n,t}$ are 
$\vert F_{1,t}\vert=t+1,\vert F_{2,t}\vert=2t+1,\vert F_{3,t}\vert=3t+2,\cdots,\vert F_{n,t}\vert=F_{n}+tF_{n+1}$.
As can be seen from the operations of determinant, sequences $\lbrace 0,-1,-1,\cdots,-F_{n-1},\cdots\rbrace$, $\lbrace 1,1,2,\cdots,F_{n},\cdots\rbrace$, $\lbrace 2,3,5, \cdots,F_{n+2},\cdots\rbrace$ and $\lbrace 3,5,8, \cdots,F_{n+3},\cdots\rbrace$ for respectively $t=-1,0,1$ and $2$.\\
Fifth, Fibonacci-Hessenberg matrix is matrix $G_{n,t}$ is accessed by replacing which element is $2$ in the $\left(1,1\right)$ of the matrix $F_{n,t}$ with $1$. The represantation of matrix $G_{n,t}$ and generalized determinants have given in the following.
\begin{equation}
G_{n,t}=\begin{pmatrix}
1 & -1 & \cdots & 0 & 0 \\
-1 & 2 & \cdots & 0 & 0 \\
\vdots & \vdots& \ddots & 0 & 0 \\
-1 & 1 & \cdots & 2 & -1 \\
1 & -1 & \cdots & -1 & t+1 
\end{pmatrix}
\end{equation}
$\vert G_{1,t}\vert=t+1,\vert G_{2,t}\vert=2t+1,\vert G_{3,t}\vert=3t+2,\cdots,\vert G_{n,t}\vert=F_{n-2}+tF_{n-1}$
If $t=-1,0,1$ and $2$ are taken, sequences $\lbrace 0,-1,-1,\cdots,-F_{n-3},\cdots \rbrace$,$\lbrace 1,1,2,\cdots,F_{n-2},\cdots\rbrace$, $\lbrace 2,3,5,\cdots,F_{n},\cdots\rbrace$, $\lbrace 3,5,8,\cdots,F_{n+1},\cdots\rbrace$ are attained for respectively.\\
Sixth, Fibonacci-Hessenberg matrix is matrix $H_{n,t}$ is accomplished by replacing which element are $-1$ in the upper diagonal of the matrix $G_{n,t}$ with $1$. The matrix $H_{n,t}$ and its generalized determinants can be given as follows:
\begin{equation}
H_{n,t}=\begin{pmatrix}
1 & 1 & \cdots & 0 & 0 \\
-1 & 2 & \cdots & 0 & 0 \\
\vdots & \vdots& \ddots & 0 & 0 \\
-1 & 1 & \cdots & 2 & 1 \\
1 & -1 & \cdots & -1 & t+1 
\end{pmatrix}
\end{equation}
$\vert H_{1,t}\vert=t+1,\vert H_{2,t}\vert=2t+3,\vert H_{3,t}\vert=5t+8,\cdots,\vert H_{n,t}\vert=F_{2n-1}+tF_{2n-2}$
The following can be presented about the sequences attained with the determiances of the matrix $H_{n,t}$: sequences $\lbrace 0,1,3,\cdots,F_{2n-3},\cdots\rbrace$, $\lbrace 1,3,8,\cdots,F_{2n-1},\cdots\rbrace$, $\lbrace 2,5,13,\cdots,F_{2n},\cdots\rbrace$ and $\lbrace 3,7,18,\cdots,L_{2n-1},\cdots\rbrace$ for respectively $t=-1,0,1$ and $2$.\\
Seventh, Fibonacci-Hessenberg matrix is matrix $K_{n,t}$ is realized by replacing which element is $1$ in the upper diagonal of the matrix $H_{n,t}$ with $2$. The matrix $K_{n,t}$  can be visualized and its generalized determinants can be given as follows:
\begin{equation}
K_{n,t}=\begin{pmatrix}
2 & 1 & \cdots & 0 & 0 \\
-1 & 2 & \cdots & 0 & 0 \\
\vdots & \vdots& \ddots & 0 & 0 \\
-1 & 1 & \cdots & 2 & 1 \\
1 & -1 & \cdots & -1 & t+1 
\end{pmatrix}
\end{equation}
$\vert K_{1,t}\vert=t+1,\vert K_{2,t}\vert=2t+3,\vert K_(,{3,t}\vert=5t+8,\cdots,\vert K_{n,t}\vert=F_{2n}+tF_{2n-}$
The sequences derived from the evaluation of the determinants of the matrix $K_{n,t}$ can be summarized as below: sequences $\lbrace 0,1,3,\cdots,F_{2n-2},\cdots\rbrace$,$\lbrace 1,3,8,\cdots,F_{2n},\cdots\rbrace$, $\lbrace 2,5,13,\cdots,F_{2n+1},\cdots\rbrace$ and $\lbrace 3,7,18,\cdots,L_{2n},\cdots\rbrace$ for respectively $t=-1,0,1$ and $2$.\\
The transitions from one to another in the Fibonacci-Hessenberg matrices given just above can be visualized as follows. In the next matrix, the changes with respect to the first matrix are depicted in red.\\
\begin{center}
$
C_{n,t}=\begin{pmatrix}
2 & 1 & \cdots & 0 & 0 \\
1 & 2 & \cdots & 0 & 0 \\
\vdots & \vdots& \ddots & 0 & 0 \\
1 & 1 & \cdots & 2 & 1 \\
1 & 1 & \cdots & 1 & t+1 
\end{pmatrix}
\longrightarrow
D_{n,t}=\begin{pmatrix}
2 & \color{red}{-1} & \cdots & 0 & 0 \\
1 & 2 & \cdots & 0 & 0 \\
\vdots & \vdots& \ddots & 0 & 0 \\
1 & 1 & \cdots & 2 & \color{red}{-1} \\
1 & 1 & \cdots & 1 & t+1 
\end{pmatrix}
\longrightarrow
$
\\ 
$
E_{n,t}=\begin{pmatrix}
\color{red}{1} & -1 & \cdots & 0 & 0 \\
1 & 2 & \cdots & 0 & 0 \\
\vdots & \vdots& \ddots & 0 & 0 \\
1 & 1 & \cdots & 2 & -1 \\
1 & 1 & \cdots & 1 & t+1 
\end{pmatrix}
\longrightarrow
F_{n,t}=\begin{pmatrix}
2 & -1 & \cdots & 0 & 0 \\
\color{red}{-1} & 2 & \cdots & 0 & 0 \\
\vdots & \vdots& \ddots & 0 & 0 \\
\color{red}{-1} & 1 & \cdots & 2 & -1 \\
1 & \color{red}{-1} & \cdots & \color{red}{-1} & t+1 
\end{pmatrix}
\longrightarrow
$
\\
$
G_{n,t}=\begin{pmatrix}
\color{red}{1} & -1 & \cdots & 0 & 0 \\
-1 & 2 & \cdots & 0 & 0 \\
\vdots & \vdots& \ddots & 0 & 0 \\
-1 & 1 & \cdots & 2 & -1 \\
1 & -1 & \cdots & -1 & t+1 
\end{pmatrix}
\longrightarrow
H_{n,t}=\begin{pmatrix}
1 & \color{red}{1} & \cdots & 0 & 0 \\
-1 & 2 & \cdots & 0 & 0 \\
\vdots & \vdots& \ddots & 0 & 0 \\
-1 & 1 & \cdots & 2 & \color{red}{1} \\
1 & -1 & \cdots & -1 & t+1 
\end{pmatrix}
\longrightarrow
$
\\
$
K_{n,t}=\begin{pmatrix}
\color{red}{2} & 1 & \cdots & 0 & 0 \\
-1 & 2 & \cdots & 0 & 0 \\
\vdots & \vdots& \ddots & 0 & 0 \\
-1 & 1 & \cdots & 2 & 1 \\
1 & -1 & \cdots & -1 & t+1 
\end{pmatrix}
$
\end{center}
The sequences generated with the determinants of the Fibonacci Hessenberg matrices mentioned previously and the recurrence relations are listed in the below table 1.
\begin{table}[h]
\caption{The sequences and recurrence relations created with the determinants of the Fibonacci Hessenberg matrices}\label{tab2}
\begin{tabular*}{\textwidth}{@{\extracolsep\fill}lccccc}
\toprule%

Matrix $C_{n,t}$ & $n=1$ & $n=2$ & $n=3$ & $\cdots$ & $n$  \\
\midrule
$t=-1$  & $0$ & $-1$ & $-1$ & $\cdots$  & $-F_{n-1}$ \\
$t=0$  & $1$ & $1$  & $2$  & $\cdots$ & $F_{n}$ \\
$t=1$  & $2$ & $3$ & $5$ & $\cdots$ & $F_{n+2}$ \\
$t=2$  & $3$ & $5$  & $8$  &   $\cdots$  & $F_{n+3}$ \\
$t$  & $t+1$ & $2t+1$ & $3t+2$ & $\cdots$ & $F_{n+1} t+F_{n}$\\
\botrule
\end{tabular*}
\begin{tabular*}{\textwidth}{@{\extracolsep\fill}lccccc}
\toprule%
Matrix $D_{n,t}$ & $n=1$ & $n=2$ & $n=3$ & $\cdots$ & $n$  \\
\midrule
$t=-1$  & $0$	& $1$ &	$2$ &	$\cdots$ &	$F_{2n-2}$ \\
$t=0$  & $1$	& $3$	 & $5$ & $\cdots$ &	$F_{2n}$ \\
$t=1$  & $2$ & $5$	& $8$	& $\cdots$ &	$F_{2n+1}$ \\
$t=2$  & $3$	& $7$ &	$11$ &	$\cdots$ &	$L_{2n}$ \\
$t$  & $t+1$ & $2t+3$ &	$3t+5$ & $\cdots$ &	$F_{2n-1} t+F_{2n}$
\\
\botrule
\end{tabular*}
\begin{tabular*}{\textwidth}{@{\extracolsep\fill}lccccc}
\toprule%
Matrix $E_{n,t}$ & $n=1$ & $n=2$ & $n=3$ & $\cdots$ & $n$  \\
\midrule
$t=-1$  & $0$	 & $1$ & $2$ &	$\cdots$ &	$F_{2n-3}$ \\
$t=0$  & $1$ &	$3$	& $5$ &	$\cdots$ &	$F_{2n-1}$ \\
$t=1$  & $2$ &	$5$ &	$8$	& $\cdots$	& $F_{2n}$ \\
$t=2$  & $3$	& $7$ &	$11$ &	$\cdots$ &	$L_{2n-1}$ \\
$t$  & $t+1$ &	$2t+3$ &	$3t+5$ & $\cdots$ &	$F_{2n-2} t+F_{2n-1}$
\\
\botrule
\end{tabular*}
\begin{tabular*}{\textwidth}{@{\extracolsep\fill}lccccc}
\toprule%
Matrix $F_{n,t}$ & $n=1$ & $n=2$ & $n=3$ & $\cdots$ & $n$  \\
\midrule
$t=-1$  & $0$ &	$-1$ &	$-1$ & $\cdots$ &	$-F_{n-1}$ \\
$t=0$  & $1$ &	$1$ &	$2$ & $\cdots$ &	$F_{n}$ \\
$t=1$  & $2$ &	$3$ &	$5$  &	$\cdots$ &	$F_{n+2}$ \\
$t=2$  & $3$ &	$5$ &	$8$	& $\cdots$ &	$F_{n+3}$ \\
$t$  & $t+1$ & $2t+1$ &	$3t+2$ & $\cdots $ &	$F_{n}+tF_{n+1}$
\\
\botrule
\end{tabular*}
\begin{tabular*}{\textwidth}{@{\extracolsep\fill}lccccc}
\toprule%
Matrix $G_{n,t}$ & $n=1$ & $n=2$ & $n=3$ & $\cdots$ & $n$  \\
\midrule
$t=-1$  & $0$ &	$-1$ &	$-1$ &	$\cdots$ &	$-F_{n-3}$ \\
$t=0$  & $1$ &	$1$ &	$2$ &	$\cdots$	& $F_{n-2}$ \\
$t=1$  & $2$ &	$3$ & $5$ &	$\cdots$ &	$F_{n}$ \\
$t=2$  & $3$ &	$5$ & $8$ &	$\cdots$ &	$F_{n+1}$ \\
$t$  & $t+1$ & $2t+1$ &	$3t+2$ & $\cdots$ &	$F_{n-2}+tF_{n-1}$
\\
\botrule
\end{tabular*}
\begin{tabular*}{\textwidth}{@{\extracolsep\fill}lccccc}
\toprule%
Matrix $H_{n,t}$ & $n=1$ & $n=2$ & $n=3$ & $\cdots$ & $n$  \\
\midrule
$t=-1$  & $0$	 & $1$ & $3$	& $\cdots$ &	$F_{2n-3}$ \\
$t=0$  & $1$ &	$3$ &	$8$ & $\cdots$ &	$F_{2n-1}$ \\
$t=1$  & $2$ &	$5$ &	$13$ &	$\cdots$ &	$F_{2n}$ \\
$t=2$  & $3$ &	$7$ &	$18$ &	$\cdots$	& $L_{2n-1}$ \\
$t$  & $t+1$ & $2t+3$ &	$5t+8$ & $\cdots$ &	$F_{2n-1}+tF_{2n-2}$
\\
\botrule
\end{tabular*}
\begin{tabular*}{\textwidth}{@{\extracolsep\fill}lccccc}
\toprule%
Matrix $K_{n,t}$ & $n=1$ & $n=2$ & $n=3$ & $\cdots$ & $n$  \\
\midrule
$t=-1$  & $0$ &	$1$ & $3$ &	$\cdots$ &	$F_{2n-2}$ \\
$t=0$  & $1$	& $3$ &	$8$ & $\cdots$ & $F_{2n}$ \\
$t=1$  & $2$ &	$5$	& $13$ & $\cdots$ &	$F_{2n+1}$ \\
$t=2$  & $3$ &	$7$ &	$18$ &	$\cdots$&	$L_{2n}$ \\
$t$  & $t+1$ & $2t+3$ &	$5t+8$ &	$\cdots$	& $F_{2n}+tF_{2n-1}$
\\
\botrule
\end{tabular*}
\end{table}
\subsection*{$C_{n,t}^{i}, D_{n,t}^{i}$ and $E_{n,t}^{i}$ Fibonacci-Hessenberg Matrices}
In previous studies, in Fibonacci-Hessenberg matrices $C_{n,t},D_{n,t}$ and $E_{n,t}$, by substituting the $i^{th}$ column by  column $1$, the matrices whose determinants are $t$-Fibonacci numbers are produced. These matrices are indicated by $C_{n,t}^{i}, D_{n,t}^{i}$ and $E_{n,t}^{i}$. Given the sequence $\lbrace \vert C_{n,t}^{i} \vert, n\geqslant 1 \rbrace$ for the matrix $C_{n,t}^{i}$, the recurrence relation is then derived as $\vert C_{n,t}^{i} \vert=tF_{n-i}+F_{n-i-1},n\geqslant i \geqslant 1$. Therefore, the recurrence relation that yields the sequence of determinants of the matrix $C_{n,t}$ is $\vert C_{n,t}^{i}\vert=t+\sum_{i=1}^{n}\vert C_{n,t}^{i} \vert$ . Likewise, for a matrix $D_{n,t}^{i}$, the recurrence relation which gives the sequence $\lbrace \vert D_{n,t}^{i} \vert, n\geqslant 1 \rbrace$ is acquired as $\vert D_{n,t}^{i}\vert=F_{2(n-i)+1}+tF_{2(n-i)},n\geqslant i \geqslant 1$. In a similar way, the recurrence relation $\vert D_{n,t}^{i}\vert=F_{2(n-i)+1}+tF_{2(n-i)},n\geqslant i \geqslant 1$ for the matrix $D_{n,t}^{i}$ gives the sequence $\lbrace\vert D_{n,t}^{i}\vert,n\geqslant 1\rbrace$. Consequently, it is the recurrence relation $\vert D_{n,t}\vert=t+\sum_{i=1}^{n}\vert D_{n,t}^{i}\vert$ that provides the sequence of determinants of the matrix $D_{n,t}$. In additionally, the recurrence relation to give the sequence $\lbrace \vert E_{n,t}^{i}\vert,n\geqslant 1\rbrace$ is determined as $\vert E_{n,t}^{i}\vert=tF_{n-i}+F_{n-i+1},n\geqslant i \geqslant 2$.Furthermore, the recurrence relation $2\vert E_{n,t} \vert=(t+1)+\sum_{i=1}^{n}\vert E_{n,t}^{i}\vert,n\geqslant 2$ satisfies the sequence of determinants of the matrix $E_{n,t}$.\\
The $t$-Fibonacci numbers of Fibonacci-Hessenberg matrices $C_{n,t}^{i}, D_{n,t}^{i}$ and $E_{n,t}^{i}$ can be seen in the table 2 for $i=1,2$ and $3$.
\begin{table}[h]
\caption{$t$-Fibonacci numbers of Fibonacci-Hessenberg matrices $C_{n,t}^{i}, D_{n,t}^{i}$ and $E_{n,t}^{i}$ for $i=1,2$ and $3$}\label{tab2}
\begin{tabular*}{\textwidth}{@{\extracolsep\fill}lccc}
\toprule%

& Matrix $C_{n,t}^{i}$ & Matrix $D_{n,t}^{i}$ & Matrix $E_{n,t}^{i}$  \\
\midrule
$i=1$ & $tF_{n-1}+F_{n-2}$ &	$F_{2n-1}+tF_{2n-2}$	& $tF_{n-1}+F_{n}$  \\
$i=2$  & $tF_{n-2}+F_{n-3}$ &	$F_{2n-3}+tF_{2n-4}$	 & $tF_{n-2}+F_{n-1}$   \\
$i=3$  & $tF_{n-3}+F_{n-4}$ &	$F_{2n-5}+tF_{2n-6}$ &	$tF_{n-3}+F_{n-2}$  \\
\botrule
\end{tabular*}
\end{table}
\subsection*{Lorentz Matrix Multiplication}
The second part highlights the key theoretical concepts which related to obtaining of Lorentz matrices . A considerable amount of literature has been published on Lorentz matrices. These studies have been widely investigated by researchers.\\
Suppose that the set of matrices with $m$ rows and $n$ columns consisting of real numbers is denoted by $R_{n}^{m}$. The matrix set $R_{n}^{m}$ is a real vector space with respect to addition and scalar multiplication. [$\color{blue}{21}$] can be referenced for Lorentz matrix multiplication. Now, suppose there are two vectors $x=(x_{1},x_{2},\cdots,x_{n})$ and $y=(y_{1},y_{2},\cdots,y_{n})\in R_{n}^{n}$. The inner product between the two vectors is described as given below:
\begin{equation}
\langle x,y \rangle=-x_{1}y_{1}+\sum_{i=2}^{n}x_{i}y_{i} 
\end{equation}
The set of the Lorentz matrices to apply the Lorentz inner product is designated by $L_{n}^{n}$.
\subsection*{Lorentz Matrix Multiplication and its Properties}
Assume that $A=[a_{ij}]\in R_{n}^{m}$ and $B=[b_{jk} ]\in R_{p}^{n}$. In addition, assume that $\cdot_L$ is a Lorentz matrix multiplication between matrices $A$ and $B$. Then,
\begin{equation}
A\cdot_L B=\left[-a_{i1}b_{1k}+\sum_{j=2}^{n}a_{ij}b_{jk}\right] 
\end{equation}
is determined. That is, $A\cdot_L B$ is a $m\times p$ type matrix. More over, every element in $A\cdot_L B$ is a Lorentz inner product. Also let $A_{i}$  be the $i^{th}$ row in $A$ and $B^{j}$ the $j^{th}$ column in $B$. Thus, the $(i,j)^{th}$ element of $A\cdot_L B$ is $\langle A_{i},B^{j}\rangle_{L}$. In this way, $A\cdot_{L}B$ can be identified as the below:
\begin{equation}
A\cdot_L B=\begin{pmatrix}
\langle A_{1},B^{1}\rangle_{L} &  \cdots & \langle A_{1},B^{j}\rangle_{L}  \\
\vdots & \ddots & \vdots \\
\langle A_{j},B^{1}\rangle_{L} & \cdots &\langle A_{j},B^{j}\rangle_{L}  
\end{pmatrix}
\end{equation}
\begin{theorem} $[\color{blue}{22}]$ The below properties related to the Lorentz matrix product are ensured:
\item[i.] $\forall A\in L_{n}^{m}, B\in L_{p}^{n}$ for $A\cdot_L (B\cdot_L C)=(A\cdot_L B)\cdot_L C$
\item[ii.] $\forall A\in L_{n}^{m}, B,C\in L_{p}^{n}$ for $A\cdot_L (B+C)=A\cdot_L B+A\cdot_LC$
\item[iii.] $\forall A,B\in L_{n}^{m}, C\in L_{p}^{n}$ for $(A+B)\cdot_L C=A\cdot_L C+B\cdot_L C$
\item[iv.] $\forall k\in R, A\in L_{n}^{m}, B\in L_{p}^{n}$ for $k(A\cdot_L B)=(kA)\cdot_L B=A\cdot_L (kB)$
\end{theorem}
With respect to Lorentz matrix multiplication, assuming that the unit matrix of type $n\times n$ is labeled $I_{n}^{n}$. That is,
\begin{equation}
I_{n}^{n}=\begin{pmatrix}
1 & 0 & \cdots & 0 & 0 \\
0 & 1 & \cdots & 0 & 0 \\
\vdots & \vdots& \ddots & 0 & 0 \\
0 & 0 & \cdots & 1 & 0 \\
0 & 0 & \cdots & 0 & 1 
\end{pmatrix}\cdot_L\begin{pmatrix}
1 & 0 & \cdots & 0 & 0 \\
0 & 1 & \cdots & 0 & 0 \\
\vdots & \vdots& \ddots & 0 & 0 \\
0 & 0 & \cdots & 1 & 0 \\
0 & 0 & \cdots & 0 & 1 
\end{pmatrix}=\begin{pmatrix}
-1 & 0 & \cdots & 0 & 0 \\
0 & 1 & \cdots & 0 & 0 \\
\vdots & \vdots& \ddots & 0 & 0 \\
0 & 0 & \cdots & 1 & 0 \\
0 & 0 & \cdots & 0 & 1 
\end{pmatrix}
\end{equation}
\begin{definition} $[\color{blue}{22}]$ $A$ is called invertible and denoted by $A^{-1}$ if for $A\in L_{n}^{n}$ if there exists a matrix $B\in L_{n}^{n}$ of type $n\times n$ satisfying the equation $A\cdot_L B=B\cdot_L A=I_{n}^{n}$.
\end{definition}
\begin{definition} $[\color{blue}{22}]$ The transpose of the matrix $A\in L_{n}^{m}$ is represented by $A^{T}$ and is then given by $A^{T}=\left[a_{ji}\right]\in L_{m}^{n}$.
\end{definition}
\begin{definition} $[\color{blue}{22}]$ $A\in L_{n}^{m}$ is called $L$-ortogonal if $A^{-1}=A^{T}$.
\end{definition}
\begin{definition} $[\color{blue}{22}]$  The $L$-determinant of the matrix $A=\left[a_{ij}\right]\in R_{n}^{n}$ is identified by $\det A$ and given by $\det A=\sum_{\sigma\in S_{n}}s\left(\sigma\right)a_{s\left( 1\right)}a_{s\left( 2\right)2}\cdots a_{s\left(n\right)n}$. Note that $S_{n}$ is a set containing all permutations of the set $\left\lbrace 1,2,\cdots,n\right\rbrace $. Furthermore, $s\left(\sigma\right)$  is the sign of the permutations of $\sigma$.
\end{definition}
\begin{theorem} $[\color{blue}{22}]$  For every $A,B\in L_{n}^{n}$, $\ det\left(A\cdot_L B\right)=-\det A \cdot \det B$ is provided.
\end{theorem}
\section{Results}\label{sec2}
The methods for obtaining a matrix have varied somewhat across this research area. Different methods have been proposed to achieve a matrix. The use of Lorentz inner product is one of the most well known methods for attaining a matrix. In this part, the matrix to which Lorentz matrix multiplication is applied is specifically denoted by $A_{L}$. Suppose that for $A=\left[a_{ij} \right]\in R_{n}^{n}$ ve $B=\left[ b_{jk}\right]\in R_{n}^{n}$, $A\cdot_L B=AB_{L}$ is defined. Here, the matrix $AB_{L}$  used is called Fibonacci-Hessenberg-Lorentz matrix. By the  Theorem 2, $\det AB_{L}=-\det A\cdot\det B$ for $A,B \in L_{n}^{n}$. By applying Lorentz matrix multiplication to the $n\times n$ type matrices discussed in the literature, it allows us to reach completely new matrices. Via the determinants of the related matrices, new sequences of numbers associated with Fibonacci numbers have been included in the literature.  
\begin{conclusion}
Suppose $C_{n,t}$ and $D_{n,t}$ are Fibonacci-Hessenberg matrices and $\cdot_L$ is the Lorentz matrix product. The determinant of the matrix derived under Lorentz matrix multiplication of the $C_{n,t}$ and $D_{n,t}$ matrices for $n\geqslant1$ yields the sequence as follows: \\
\begin{center}
$
\vert CD_{n,t}\vert=\lbrace -t^{2}-2t-1,-4t^{2}-8t-3,-15t^{2}-34t-16,\cdots,-F_{n+1}F_{2n-1} t^{2}-\left(F_{n+1}F_{2n}+F_{n}F_{2n-1}\right)t-F_{n} F_{2n},\cdots\rbrace
$
\end{center}
\end{conclusion}
\begin{proof}
For $n=1$, $C_{1,t}=\begin{pmatrix} t+1 \end{pmatrix}$   and $D_{1,t}=\begin{pmatrix} t+1 \end{pmatrix}$ are provided. Then,
\begin{center}
$
C_{1,t} \cdot_L D_{1,t}=CD_{1,t}=\begin{pmatrix} -t^{2}-2t-1 \end{pmatrix}
$
\end{center}
is achieved. By Theorem 2,
\begin{center}
$
\det CD_{1,t}=-\det C_{1,t}\cdot \det D_{1,t}
=-\left(t+1\right)\cdot\left(t+1\right)
=-t^{2}-2t-1
$
\end{center}
is found. For $n=2$, $C_{2,t}=\begin{pmatrix} 2 & 1 \\1 & t+1 \end{pmatrix}$ and $D_{2,t}=\begin{pmatrix} 2 & -1 \\1 & t+1 \end{pmatrix}$ are satisfied. 
\begin{center}
$
CD_{2,t}=C_{2,t}\cdot_L D_{2,t}
=\begin{pmatrix} 2 & 1 \\1 & t+1 \end{pmatrix}\cdot_L \begin{pmatrix} 2 & -1 \\1 & t+1 \end{pmatrix}
=\begin{pmatrix} -3 & t+3 \\t-1 & t^{2}+2t+2 \end{pmatrix}
$
\end{center}
is then obtained under Lorentz matrix multiplication. The determinant of $CD_{2,t}$,
\begin{center}
$
\det CD_{2,t}=-\det C_{2,t}\cdot \det D_{2,t}
=-\left(2t+1\right)\cdot\left(2t+3\right)
=-4t^{2}-8t-3
$
\end{center}
is attained. The matrices for $n=3$ are equal to 
$C_{3,t}=\begin{pmatrix} 2 & 1 & 0 \\1 & 2 & 1 \\ 1 & 1 & t+1 \end{pmatrix}$ and $D_{3,t}=\begin{pmatrix} 2 & 1 & 0 \\1 & 2 & 1 \\ 1 & 1 & t+1 \end{pmatrix}.$
Applying the Lorentz matrix multiplication to these resulting matrices produces 
\begin{center}
$
CD_{3,t}=C_{3,t}\cdot_L D_{3,t}
=\begin{pmatrix} -3 & 4 &-1 \\ 1 & 6& t-1 \\ t & t+4 & t^{2}+2t \end{pmatrix}
$.
\end{center}
The determinant given in the Theorem 2 provides 
\begin{center}
$
\det CD_{3,t}=-\det C_{3,t}\cdot \det D_{3,t}
=-\left(3t+2\right)\cdot\left(5t+8\right)
=-15t^{2}-34t-16
$.
\end{center}
It is determined as 
$
C_{n,t}=\begin{pmatrix}
2 & 1 & \cdots & 0 & 0 \\
1 & 2 & \cdots & 0 & 0 \\
\vdots & \vdots& \ddots & 0 & 0 \\
1 & 1 & \cdots & 2 & 1 \\
1 & 1 & \cdots & 1 & t+1 
\end{pmatrix}
$
and
$
D_{n,t}=\begin{pmatrix}
2 & -1 & \cdots & 0 & 0 \\
1 & 2 & \cdots & 0 & 0 \\
\vdots & \vdots& \ddots & 0 & 0 \\
1 & 1 & \cdots & 2 & -1 \\
1 & 1 & \cdots & 1 & t+1 
\end{pmatrix}
$
for the general $n$ term.
By Equation $C_{n,t}\cdot_L D_{n,t}=CD_{n,t}$ and Theorem 2, 
\begin{center}
$
\det CD_{n,t}=-\det C_{n,t}\cdot\det D_{n,t}
=-\left(F_{n+1}t+F_{n}\right)\cdot\left(F_{2n-1}t+F_{2n}\right)
=-F_{n+1}F_{2n-1}t^{2}-\left(F_{n+1}F_{2n}+F_{n}F_{2n-1} \right)t-F_{n}F_{2n}
$
\end{center}
is derived. Consequently, the below sequences of numbers are arrived at: 
\begin{center}
$
\vert CD_{n,t}\vert=\lbrace-t^{2}-2t-1,-4t^{2}-8t-3,-15t^{2}-34t-16,\cdots,-F_{n+1} F_{2n-1}t^{2}-\left(F_{n+1}F_{2n}+F_{n}F_{2n-1}\right)t-F_{n}F_{2n},\cdots\rbrace
$.
\end{center}
\end{proof}
\begin{conclusion}
Let $E_{n,t}$ and $F_{n,t}$ be the Fibonacci-Hessenberg matrices and $\cdot_L$ be the Lorentz matrix multiplication. For $n\geqslant1$, the determinant of the matrix $EF_{n,t}$ generated under the Lorentz matrix multiplication of $E_{n,t}$ and $F_{n,t}$ gives the next sequence: \\
\begin{center}
$
\vert EF_{n,t}\vert=\lbrace -t^{2}-2t-1,-4t^{2}-8t-3,-15t^{2}-34t-16,\cdots,-F_{2n-2} F_{n+1} t^{2}-\left(F_{2n-2}F_{n}+F_{n+1}F_{2n-1}\right)-F_{2n-1}F_{n},\cdots\cdots\rbrace
$
\end{center}
\end{conclusion}
\begin{proof}
Note that for $n=1$, $E_{1,t}=\begin{pmatrix} t+1 \end{pmatrix}$ and $F_{1,t}=\begin{pmatrix} t+1 \end{pmatrix}$ are achieved. After that,
\begin{center}
$
E_{1,t} \cdot_L F_{1,t}=EF_{1,t}=\begin{pmatrix} -t^{2}-2t-1 \end{pmatrix}
$
\end{center}
is obtained. Hence, applying the determinant operation according to the Lorentz matrix multiplication, it is clear that 
\begin{center}
$
\det EF_{1,t}=-\det E_{1,t}\cdot \det F_{1,t}
=-\left(t+1\right)\cdot\left(t+1\right)
=-t^{2}-2t-1
$.
\end{center}
In the following, related Lorentz matrix and its determinant value are presented for matrices 
$E_{2,t}=\begin{pmatrix} 2 & -1 \\ 1 & t+1 \end{pmatrix}$ and $F_{2,t}=\begin{pmatrix} 2 & -1 \\ -1 & t+1 \end{pmatrix}$:
\begin{center}
$E_{2,t}\cdot_L F_{2,t}=EF_{2,t}=\begin{pmatrix} -3 & -t+1 \\ -t-3 & t^{2}+2t+2 \end{pmatrix}$
\end{center}
\begin{center}
$\det EF_{2,t}=-\det E_{2,t}\cdot \det F_{2,t}
=-4t^{2}-8t-3$
\end{center}
By applying similar operations,
$E_{3,t}=\begin{pmatrix} 2 & 1 & 0 \\ 1 & 2 & 1 \\ 1 & 1 & t+1 \end{pmatrix}$
and 
$
F_{3,t}=\begin{pmatrix} 2 & -1 & 0 \\ 1 & 2 & -1 \\ 1 & -1 & t+1 \end{pmatrix}
$
matrices are reached for $n=3$. Subsequently, the corresponding Lorentz matrix is equal to
\begin{center}
$E_{3,t}\cdot_L F_{3,t}=EF_{3,t}=\begin{pmatrix} -3 & 0 & 1 \\ -5 & 6 & -t-3 \\ t-2 & -t+2 & t^{2}+2t \end{pmatrix}$.
\end{center}
On the other hand, the determinant for matrix $EF_{3,t}$ is 
\begin{center}
$\det EF_{3,t}=-\det E_{3,t}\cdot \det F_{3,t}
=-(5t+8)\cdot(3t+2)
=-15t^{2}-34t-16$.
\end{center}
A Lorentz matrix is also generated from the matrices 
$
E_{n,t}=\begin{pmatrix}
1 & -1 & \cdots & 0 & 0 \\
1 & 2 & \cdots & 0 & 0 \\
\vdots & \vdots& \ddots & 0 & 0 \\
1 & 1 & \cdots & 2 & -1 \\
1 & 1 & \cdots & 1 & t+1 
\end{pmatrix}
$
and
$
F_{n,t}=\begin{pmatrix}
2 & -1 & \cdots & 0 & 0 \\
-1 & 2 & \cdots & 0 & 0 \\
\vdots & \vdots& \ddots & 0 & 0 \\
-1 & 1 & \cdots & 2 & -1 \\
1 & -1 & \cdots & -1 & t+1 
\end{pmatrix}
$
of type $n\times n$ and its determinant is then assigned as
\begin{center}
$\det EF_{n,t}=-\det E_{n,t}\cdot \det F_{n,t}
=-\left(F_{2n-2}t+F_{2n-1}\right)\cdot\left(F_{n}+tF_{n+1} \right)
=-F_{2n-2}F_{n+1}t^{2}-\left(F_{2n-2}F_{n}+F_{n+1}F_{2n-1}\right)t-F_{2n-1}F_{n}$.
\end{center}
As a result, the following sequences of numbers are produced:  
\begin{center}
$ \vert EF_{n,t}\vert=\lbrace-t^{2}-2t-1,-4t^{2}-8t-3,-15t^{2}-34t-16,\cdots,-F_{2n-2} F_{n+1}t^{2}-\left( F_{2n-2}F_{n}+F_{n+1}F_{2n-1}\right)t-F_{2n-1}F_{n},\cdots \rbrace $.
\end{center}
\end{proof}
\begin{conclusion} For $n\geqslant1$, the determinant of the matrix $GH_{n,t}$ formed by the Lorentz matrix multiplication of $G_{n,t}$ and $H_{n,t}$ produces the subsequent sequence:\\
\begin{center}
$
\vert GH_{n,t}\vert=\lbrace-t^{2}-2t-1,-4t^{2}-8t-3,-15t^{2}-34t-16,\cdots,-F_{2n-2} F_{n+1}t^{2}-\left( F_{2n-2}F_{n}+F_{n+1}F_{2n-1}\right)t-F_{2n-1}F_{n},\cdots\rbrace
$
\end{center}
\end{conclusion}
\begin{proof}
The matrices $G_{n,t}$ and $H_{n,t}$ are given as $G_{1,t}=\begin{pmatrix} t+1 \end{pmatrix}$, $H_{1,t}=\begin{pmatrix} t+1 \end{pmatrix}$ and $G_{2,t}=\begin{pmatrix} 2 & -1 \\ -1 & t+1 \end{pmatrix}$, $H_{2,t}=\begin{pmatrix} 2 & 1 \\ -1 & t+1 \end{pmatrix}$ for $n=1$ and $2$ respectively. It is obvious that the matrices $GH_{1,t}=\begin{pmatrix} -t^{2}-2t-1 \end{pmatrix}$ and $GH_{2,t}=\begin{pmatrix} -3 &-t-3 \\ -t+1 & t^{2}+2t+2 \end{pmatrix}$, respectively.
By the Lorentz matrix multiplication determinant operation, their determinants 
\begin{center}
$
\det GH_{1,t}=-\det G_{1,t}\cdot \det H_{1,t}
=-\left(t+1\right)\cdot \left(t+1\right) 
=-t^{2}-2t-1
$
\end{center}
and
\begin{center}
$
\det GH_{2,t}=-\det G_{2,t}\cdot\det H_{2,t}
=-\left(2t+1\right)\cdot \left(2t+3\right)
=-4t^{2}-8t-3
$
\end{center}
are attained. For $n=3$, $G_{3,t}=\begin{pmatrix} 2 & -1 & 0 \\ -1 & 2 & -1 \\ 1 & -1 & t+1 \end{pmatrix}$ and $H_(3,t)=\begin{pmatrix} 2 & 1 & 0 \\ 1 & 2 & 1 \\ 1 & -1 & t+1\end{pmatrix}$ are determined. Here,
\begin{center}
$
G_{3,t}\cdot_L H_{3,t}=GH_{3,t}=\begin{pmatrix} -3 & -4 & -1 \\ -1 & 6 & -t+1 \\ t & -t-4 & t^{2}+2t \end{pmatrix}
$
\end{center}
is found. Similarly,
\begin{center}
$
\det GH_{3,t}=-\det G_{3,t}\cdot\det H_{3,t}
=-\left(3t+2\right)\cdot \left(5t+8 \right)
=-15t^{2}-34t-16
$
\end{center}
is realized. For general term $n$, 
$
G_{n,t}=\begin{pmatrix}
1 & -1 & \cdots & 0 & 0 \\
-1 & 2 & \cdots & 0 & 0 \\
\vdots & \vdots& \ddots & 0 & 0 \\
-1 & 1 & \cdots & 2 & -1 \\
1 & -1 & \cdots & -1 & t+1 
\end{pmatrix}
$
and
$
H_{n,t}=\begin{pmatrix}
1 & 1 & \cdots & 0 & 0 \\
-1 & 2 & \cdots & 0 & 0 \\
\vdots & \vdots& \ddots & 0 & 0 \\
-1 & 1 & \cdots & 2 & 1 \\
1 & -1 & \cdots & -1 & t+1 
\end{pmatrix}
$
considered. By equation $G_{n,t}\cdot_L H_{n,t}=GH_{n,t}$ and Theorem 2,
\begin{center}
$
\det GH_{n,t}=-\det G_{n,t} \cdot \det H_{n,t}
=-\left(F_{n-2}+tF_{n-1} \right)\cdot \left( F_{2n-1}+tF_{2n-2}\right)
=-F_{n-1}F_{2n-2}t^{2}-\left(F_{n-2}F_{2n-2}+F_{n-1} F_{2n-1}\right)t-F_{n-2}F_{2n-1}
$
\end{center}
is found. So, the following sequences of numbers are derived:
\begin{center}
$
\vert GH_{n,t}\vert =\lbrace -t^{2}-2t-1,-4t^{2}-8t-3,-15t^{2}-34t-16,\cdots,-F_{2n-2} F_{n+1}t^{2}-\left(F_{2n-2}F_{n}+F_{n+1}F_{2n-1}\right)t-F_{2n-1}F_{n},\cdots\rbrace .
$
\end{center}
\end{proof}
The table of sequences with determinants of Fibonacci-Hessenberg-Lorentz matrices $\vert CD_{n,t}\vert,\vert EF_{n,t} \vert$ and $\vert GH_{n,t}\vert$  for different values of  $n$ and $t$ can be illustrated as table 3.
\begin{table}[h]
\caption{The table related to sequences obtaining by determinants of Fibonacci-Hessenberg-Lorentz matrices $\vert CD_{n,t}\vert, \vert EF_{n,t}\vert $ and $\vert GH_{n,t}\vert$ for different values of $n$ and $t$ }\label{tab3}
\begin{tabular*}{\textwidth}{@{\extracolsep\fill}lcccccc}
\toprule%
& \multicolumn{3}{@{}c@{}}{Element} \\\cmidrule{2-6}%
$\vert CD_{n,t} \vert$ & $1^{th}$  & $2^{nd}$ & $3^{rd}$  & $\cdots$ & $ n^{th}$   \\
\midrule
$t=-1$  & $0$	 & $1$	& $3$ & $\cdots $	 &   $-F_{n+1}F_{2n-1}+\left(F_{n+1}F_{2n}+F_{n}F_{2n-1}\right)-F_{n}F_{2n}$ \\
$t=0$  & $-1$	& $-3$ &	$-16$ &	$\cdots$ &	$-F_{n}F_2{n}$  \\
$t=1$  & $-4$ &	$-15$ &	$-65$ & $\cdots$ &	$-F_{n+1} F_{2n-1}-\left(F_{n+1}F_{2n}+F_{n}F_{2n-1}\right)-F_{n} F_{2n}$  \\
$t=2$  & $-9$	& $-35$ &	$-144$ & $\cdots$ &	$-4F_{n+1} F_{2n-1}-2\left(F_{n+1}F_{2n}+F_{n}F_{2n-1}\right)-F_{n} F_{2n}$  \\
\\
\botrule
\end{tabular*}
\begin{tabular*}{\textwidth}{@{\extracolsep\fill}lcccccc}
\toprule%
& \multicolumn{3}{@{}c@{}}{Element} \\\cmidrule{2-6}%
$\vert EF_{n,t} \vert$ & $1^{th}$  & $2^{nd}$ & $3^{rd}$  & $\cdots$ & $ n^{th}$   \\
\midrule
$t=-1$  & $0$	&  $1$ &	$3$	& $\cdots$ &	$-F_{2n-2} F_{n+1}+\left( F_{2n-2}F_{n}+F_{n+1}F_{2n-1}\right)-F_{2n-1}F_{n}$ \\
$t=0$  & $-1$	& $-3$ & $-16$ & $\cdots$ &	$-F_{2n-1}F_{n}$ \\
$t=1$  & $-2$	& $-15$ &	$-65$	& $\cdots$	& $-F_{2n-2}F_{n+1}-\left(F_{2n-2}F_{n}+F_{n+1}F_{2n-1}\right)-F_{2n-1}F_{n}$  \\
$t=2$  & $-9$ &	$-35$ &	$-144$ &  $\cdots$ &	$-4F_{2n-2} F_{n+1}-2\left(F_{2n-2}F_{n}+F_{n+1}F_{2n-1}\right)-F_{2n-1}F_{n}$  \\
\\
\botrule
\end{tabular*}
\begin{tabular*}{\textwidth}{@{\extracolsep\fill}lcccccc}
\toprule%
& \multicolumn{3}{@{}c@{}}{Element} \\\cmidrule{2-6}%
$\vert GH_{n,t} \vert$ & $1^{th}$  & $2^{nd}$ & $3^{rd}$  & $\cdots$ & $ n^{th}$   \\
\midrule
$t=-1$  & $0$ &	$1$ &	$3$ & $\cdots$ &	$-F_{n-1} F_{2n-2}+\left(F_{n-2}F_{2n-2}+F_{n-1}F_{2n-1}\right)-F_{n-2}F_{2n-1}$ \\
$t=0$  & $-1$	& $-3$ & $-16$ &	$\cdots$	& $-F_{n-2} F_{2n-1}$ \\
$t=1$  & $-2$	& $-15$	& $-65$ &	$\cdots$ &	$-F_{n-1} F_{2n-2}t^{2}-\left( F_{n-2}F_{2n-2}+F_{n-1}F_{2n-1}\right)-F_{n-2}F_{2n-1}$ \\
$t=2$  & $-9$	& $-35$ &	$-144$ &	$\cdots$ &	$-F_{n-1}F_{2n-2}t^{2}-2\left(F_{n-2}F_{2n-2}+F_{n-1} F_{2n-1}\right)-F_{n-2}F_{2n-1}$  \\
\\
\botrule
\end{tabular*}
\end{table}
\section{Conclusion}\label{sec3}
In this paper, Hessenberg and Tridiagonal matrices are extensively analyzed. In besides, $\left(t,n\right)$-Fibonacci numbers derived from the determinants of Hessenberg and Tridiagonal matrices of type $n\times n$ are investigated. In contrary to other works, if $A$ and $B$ are Fibonacci-Hessenberg matrices, then $AB_{L}$ is generated by Lorentz matrix multiplication. Therefore, Fibonacci-Hessenberg-Lorentz matrices are formed for $n=1,2$ and $3$. New sequences of numbers related to the determinants of Fibonacci-Hessenberg-Lorentz matrices are determined.

\end{document}